\input amstex


\def\b1{\text{\bf 1}}

\def\CM{{\Cal M}}

\def\CV{{\Cal V}}

\def\Hom{\text{Hom}}

\def\#{\,\check{}}

\def\id{\text{id}}


\def\hra{\hookrightarrow}
\def\iso{\buildrel\sim\over\rightarrow} 

\def\lra{\longrightarrow}

\parskip=6pt

\documentstyle{amsppt}
\document
\magnification=1100
\NoBlackBoxes

\centerline {\bf Remarks on Grothendieck's standard conjectures} 

\bigskip

\centerline {A.~Beilinson}

\medskip

\centerline{ The University of Chicago}

\bigskip

We show that Grothendieck's standard conjectures (over a field of characteristic zero) follow from either of two other motivic conjectures, namely, that of existence of the motivic t-structure and (a weak version of) Suslin's  Lawson homology conjecture.
I am grateful to H.~Esnault, E.~Friedlander, and B.~Kahn for a stimulating exchange of letters.
\bigskip

\centerline{\bf \S 1. The motivic t-structure conjecture yields standard conjectures}

{\bf 1.1.}  Below $k$ is our base field, $\CV ar_k$ is the category of smooth varieties over  $k$, and  $DM_k$ is the triangulated category of geometric motives with $\Bbb Q$-coefficients over $k$,\footnote{For Voevodsky's construction of $DM_k$,  see \cite{A} (an introduction) and \cite{V}, \cite{MVW}, \cite{D\'eg},  \cite{BV} (detailed expositions). There are also equivalent approaches of M.~Hanamura and M.~Levine.} so we have the motive functor $M : \CV ar_k \to DM_k$. Recall that $DM_k$ is an idempotently complete triangulated rigid tensor $\Bbb Q$-category, and $M$ yields a fully faithful  embedding of tensor $\Bbb Q$-categories $CHM_k \hra DM_k$, were $CHM_k$ is the category of Chow motives with $\Bbb Q$-coefficients over $k$. 

Denote by $\CV ec_{\Bbb Q}$,  $\CV ec_{\Bbb Q}$  the categories of finite-dimensional $\Bbb Q$- and $\Bbb Q_\ell$-vector spaces. 
For $\ell$ prime to the characteristic of $k$ one has the $\ell$-adic realization functors $r_{\Bbb Q_\ell}: DM_k \to D^b (\CV ec_{\Bbb Q_\ell })$. For $k$ of characteristic 0,  each embedding $\iota : k\hra \Bbb C$ yields the  Betti realization functor $r_\iota : DM_k \to D^b (\CV ec_{\Bbb Q})$. These are tensor triangulated functors; there are canonical identifications $r_\iota \otimes \Bbb Q_\ell \iso r_{\Bbb Q_\ell}$.

Let $r$ be one of the realization functors. For a variety $X$, the vector spaces $H^a_r (X):=H^a r( M(X)) $ are appropriate homology of $X$; they vanish unless $-2d_X \le a \le 0$ (here $d_X :=\dim X$).

{\bf 1.2.} Let $\mu$ be a t-structure  on $DM_k$. Denote by $ DM^{\le 0}_k$, $ DM^{\ge 0}_k$ the positive and negative parts of $DM_k$, by $\CM_k$
   its heart,  and by  ${}^\mu \! H  :  DM_k \to \CM_k$  the cohomology functor.

{\bf Definition.} {\it $\mu$  is said to be motivic if it is non-degenerate\footnote{Which means that the cohomology functor ${}^\mu\! H^\cdot$ is conservative, i.e., a morphism $f$ in $DM_k$ is an isomorphism if   all ${}^\mu\!  H^a (f)$ are isomorphisms. Since we are in the triangulated setting, this amounts to the property that  any object $P$ with
${}^\mu\! H^a (P)=0$ for all $a$ equals $0$. } 
 and compatible with $\otimes$ and $r$ (i.e., $\otimes$ and $r$ are t-exact).}

{\bf Conjecture} (cf.~\cite{A} Ch.~21). {\it A motivic t-structure exists. }

Assuming the conjecture, let us deduce some of its corollaries.

{\bf 1.3.} For a motivic $\mu$, let  ${}^\mu r : \CM_k \to \CV ec_{\Bbb Q_{(\ell )}}$ be the restriction of $r$ to $\CM_k$.

{\bf Observation.} {\it $\CM_k$  is a Tannakian $\Bbb Q$-category, and ${}^\mu r$ is a fiber functor.}

{\it Proof.} $\CM_k$ is an abelian tensor $\Bbb Q$-category; the endomorphism ring of its unit object $\Bbb Q (0)$ is $\Bbb Q$. It is rigid (for t-exactness of $\otimes$ implies that the duality is t-exact). Since ${}^\mu r$ is an exact tensor functor,  we are done (see e.g.~\cite{Del2} 2.8).  \hfill$\square$

{\bf Corollary.} {\it ${}^\mu r$ is faithful (hence conservative), and every object of $\CM_k$ has finite length.
 The functor $r$ is conservative.}

{\it Proof.} The first assertion is a part of the Tannakian story. It implies the second one, for our t-structure is non-degenerate and  ${}^\mu r\, {}^\mu \! H^\cdot = H^\cdot r$ (since $r$ is t-exact).
\hfill$\square$

{\bf Corollary.} {\it  (i) Any object $P$ of $ DM_k$ has only finitely many non-zero cohomology objects ${}^\mu\! H^a P$ (i.e., $\mu$ is bounded). \newline (ii) $P$ lies in 
 $DM_k^{\le 0}$, resp.~$ DM_k^{\ge 0}$, if and only if the complex $r(P)$ has trivial positive, resp.~negative, cohomology.  }
 
 {\it Proof.} (i) Since $r$ is t-exact, one has $H^\cdot r(P)={}^\mu r\, {}^\mu\! H (P)$. The first assertion follows then from the conservativity of ${}^\mu r$. (ii)
 Since $\mu$ is non-degenerate, $P$ lies in 
 $DM_k^{\le 0}$, resp.~$ DM_k^{\ge 0}$, if and only if it has trivial positive, resp.~negative, cohomology ${}^\mu\! H^\cdot P$. We are done by  the conservativity of ${}^\mu r$.  \hfill$\square$

{\it Remarks.} (a) By (ii) above, a motivic t-structure is unique (for given $r$). If char $k=0$, then the choice of $r$ is irrelevant (indeed, since $r_{\Bbb Q_\ell }= r_\iota \otimes \Bbb Q_\ell$, the conditions from (ii) do not depend on the choice of $r$). So we can call $\mu$ {\it the} motivic t-structure. 
(b) By (ii) above, the Tate motive $\Bbb Q (1)$  lies in $ \CM_k$, so the Tate twist is t-exact.

{\bf 1.4.} Let $X$ be a smooth projective variety. Then $CH^n (X)_{\Bbb Q}=\Hom (M(X), \Bbb Q (n)[2n])$, and  the intersection product on $CH^\cdot (X)_{\Bbb Q}$ comes from the canonical coalgebra structure on $M(X)$    (the coproduct is the diagonal map  $M(X)\to M(X\times X)=M(X)\otimes M(X)$) and the evident algebra structure on $\Bbb Q (\cdot )[2\cdot ]$. 
Thus the Chow ring acts on $M(X)(\cdot )[2\cdot ]$; explicitly, the multiplication by $c\in CH^n (X)_{\Bbb Q}$ is the composition $ \cap c$ of  $M(X)\to M(X)\otimes M(X)\buildrel{\id_{M(X)}\otimes c}\over\lra M(X)(n)[2n]$.

 Let  $L\in CH^1 (X)_{\Bbb Q}=\Hom (M(X),\Bbb Q (1)[2])$ be the class of hyperplane section, so we have the morphisms $\cap L^i : M(X)\to M(X)(i)[2i]$, $i\ge 0$.

{\bf Proposition.} {\it (i) (hard Lefschetz) For any $i\ge 0$ the morphism \newline $\cap L^i : {}^\mu\! H^{-i-d_X} M(X)\to {}^\mu\!H^{i-d_X}M(X)(i)$ is an isomorphism. \newline
(ii) The object $M(X)$  is isomorphic to the direct sum of its cohomology objects: $M(X)\simeq \oplus\, {}^\mu\! H^a M(X)[-a]$.
\newline
(iii) (primitive  decomposition) There is a unique collection of subobjects ${}^\mu\! P^a (X) \subset {}^\mu \! H^{a}M(X)$,  $-2d_X\le a\le -d_X$, such that for every $b$ the maps $\cap L^i$ provide isomorphisms $\mathop\oplus\limits_{b/2 +d_X \ge i\ge \text{max} ( b+d_X, 0)} {}^\mu \! P^{b-2i} (-i)\iso {}^\mu\! H^b M(X)$. }

{\it Proof.} Applying ${}^\mu r$ to the morphism in (i), we get an isomorphism from the usual hard Lefschetz theorem for $H^\cdot_r$, and (i) follows since ${}^\mu r$ is conservative; (i) implies (ii) by \cite{Del3}; (iii) follows from (i) by a usual linear algebra argument. \hfill$\square$

{\it Remarks.} (a)   ${}^\mu r$ sends the decomposition of (iii) to a similar decomposition of ${}^\mu r\, {}^\mu\! H^\cdot M(X)=  H^\cdot r(M(X)) $, which is the
usual primitive decomposition of the homology $H^\cdot_r (X)$. 
(b) The decomposition in (ii) is usually non unique.

{\bf Corollary.} {\it The standard conjectures of Lefschetz and K\"unneth type\footnote{See e.g.~\cite{A} Ch.~5.} are true.}

{\it Proof.} Recall that $\Hom (M(X),M(X)(a)[2a])=CH^{d_X +a}(X\times X)_{\Bbb Q}$, and for any $\lambda : M(X)\to M(X)(a)[2a]$ the map 
${}^\mu r (\lambda )$ is the action on $H^\cdot_r (X)$ of the corresponding algebraic correspondence. The K\"unneth type assertion follows if we take for $\lambda$ the projector to any of the components of the decomposition from (ii). To deduce the Lefschetz type assertion, consider $\lambda : M(X)\to M(X)(-i)[-2i]$ whose only non-zero component with respect to the decomposition in (ii) is the inverse to the isomorphism from (i). \hfill$\square$

{\bf 1.5.} From now on we assume that $k$ has characteristic 0. Then (see
  \cite{A} Ch.~5) the standard conjecture of Lefschetz type implies all the standard conjectures.

{\bf Proposition.} {\it For $X$ as in 1.4, the objects ${}^\mu\! H^a M(X)$ of $\CM_k$
are semi-simple. }

{\it Proof} (cf.~\cite{J}). Due to the primitive decomposition (see 1.4(iii)) it suffices to check that the objects ${}^\mu\! P^a (X)$, $-2d_X\le a\le -d_X$, are semi-simple, which means that any subobject $Q \subset {}^\mu\! P^a (X)$ admits a complement $Q^\perp$. 

Recall that the dual  $M(X)^*$ identifies naturally with $M(X)(-d_X)[-2d_X]$. The corresponding canonical  pairing $ M(X) \otimes M(X) \to \Bbb Q (d_X)[2d_X]$ yields a pairing $(\cdot ,\cdot ):{}^\mu\! H^a  M(X)\otimes {}^\mu\! H^{-a-2d_X}  M(X)\to \Bbb Q (d_X ) $. Consider the pairing $(\cdot ,\cdot )_L := (\cdot ,\cap L^{-a-d_X} \cdot ) : {}^\mu\! P^a (X)\otimes {}^\mu\! P^a (X)\to \Bbb Q (2a)$. Let $Q^\perp \subset {}^\mu\! P^a (X)$ be the orthogonal complement to $Q$ for this pairing. 

It remains to show that $Q\oplus Q^\perp \iso {}^\mu\! P^a (X)$. It suffices to check this after applying a fiber functor ${}^\mu r_\iota$. Then $(\cdot ,\cdot )_L$ becomes the usual polarization pairing on primitive cycles, and ${}^\mu r_\iota (Q^\perp )$ is the orthogonal complement to ${}^\mu r_\iota (Q)$ with respect to the polarization. Since $(\cdot ,\cdot )_L$  is non-degenerate on ${}^\mu r_\iota (Q)$ (as on every Hodge substructure) by the Hodge index theorem,  we are done.  \hfill$\square$

 {\bf Corollary.} {\it Each irreducible object of $\CM_k$ can be realized as a Tate twist of a direct summand of some ${}^\mu\! H^{-d_X} M(X)$ where $X$ is projective and smooth.} 
 
 {\it Proof.} Every irreducible object can be realized as a subquotient of a Tate twist of some $ {}^\mu\! H^{a} M(Y)$, $Y\in \CV ar_k$. Writing $ {}^\mu\! H^{\cdot} M(Y)$ in terms of cohomology of smooth projective varieties \cite{Del1}, we see that it can be realized as a subquotient of a Tate twist of some $ {}^\mu\! H^{a} M(Y)$ with $Y$ projective and smooth. By Lefschetz, one can realize it as a subquotient of a Tate twist of  $ {}^\mu\! H^{-d_X} M(X)$ with $X$ projective and smooth. We are done by the proposition. \hfill$\square$

{\bf 1.6.} The motivic t-structure can be characterized in purely geometric terms (without the reference to $r$). Namely, consider a filtration $DM_{k (0)}\subset DM_{k (1)} \subset \ldots $ on $DM_k$ where $DM_{k (n)}$ is the thick subcategory of $DM_k$ generated by all motives of type $M(X)(a)$ with $d_X \le n$. 

{\bf Proposition.} {\it The motivic t-structure is a unique t-structure compatible with the filtration $DM_{k (\cdot )}$ and such that the heart of the induced t-structure on any successive quotient $DM_{k (n)}/DM_{k (n-1)}$ 
contains all objects $M(X)[-n]$, $X$ is projective smooth of dimension $n$.}

{\it Proof.}
By 1.4(ii),  1.5, and the argument from \cite{Del1}, $DM_{k (n)}$ is generated by all irreducible 
objects of $\CM_k$ that can be realized as a Tate twist of a direct summand of some ${}^\mu\! H^{-d_X} M(X)$ where $X$ is projective and smooth of dimension $\le n$. This shows that $\mu$ induces a t-structure on each  $DM_{k (n)}$, i.e., $\mu$ is compatible with the filtration $DM_{k (\cdot )}$. 
In such a situation, the t-structure on $DM_k = \cup DM_{k(n)}$ is uniquely determined by the t-structures induced on the successive quotients  $DM_{k (n)}/DM_{k (n-1)}$. 

The heart of the t-structure on $DM_{k (n)}$ is $\CM_k \cap DM_{k (n)}$, which is  the Serre subcategory of $\CM_k$ generated by irreducibles  occuring as direct summands of some ${}^\mu\! H^{-d_X} M(X)(a)$ where $X$ is projective and smooth of dimension $\le n$.  The irreducibles in the heart of 
$DM_{k (n)}/DM_{k (n-1)}$ are the images of those of them with $d_X=n$. For such an $X$ the image of $M(X)(a)[-n]$ in $DM_{k (n)}/DM_{k (n-1)}$ equals the image of ${}^\mu\! H^{-n} M(X)(a)$ (since ${}^\mu\! H^{\neq -n} M(X)(a)\in DM_{k (n-1)}$). Since the  t-structure on $DM_{k (n)}/DM_{k (n-1)}$
 is bounded and its heart is Artinian, it is uniquely defined by the datum of irreducible objects in its heart,\footnote{ Precisely, $(DM_{k (n)}/DM_{k (n-1)})^{\le 0}$ is the left orthogonal complement to the set of objects $M(X)(a)[\ell ]$, $X$ is projective smooth of dimension $n$, $\ell <-n$; $(DM_{k (n)}/DM_{k (n-1)})^{\ge 0}$ is the right orthogonal complement to the set of $M(X)(a)[\ell ]$, $X$ is projective smooth of dimension $n$, $\ell >-n$.} q.e.d.
  \hfill$\square$

{\bf 1.7. Proposition.} {\it (i) Objects of $\CM_k$ carry a natural finite increasing filtration $W_\cdot$ such that each morphism is strictly compatible with $W$. It is characterized by the next property: an irreducible object $P$ has weight $m$, i.e., has property $W_m P=P$, $W_{m-1}P=0$, if and only if it occurs in some ${}^\mu\! H^i M(X)(a)$ where $X$ is smooth projective and $m=i-2a$. 
\newline
(ii) If $P$, $Q$ are irreducible of weights $m$, $n$, then $\Hom (P,Q[\ell ])=0$ for $\ell >m-n$.}

{\it Proof} (cf.~\cite{Del4} 3.8). It suffices to check that if irreducible objects $P$, $Q$ occur in, respectively,  ${}^\mu\! H^i M(X)(a)$ and ${}^\mu\! H^j M(Y)(b)$, $X$ and $Y$ are smooth projective, then $\Hom (P,Q[\ell ])=0$ for $\ell >(i-2a)-(j-2b)$. By Lefschetz, we can assume that $i=-d_X$, $j=-d_Y$. By 1.5 and 1.4(ii), $\Hom (P,Q[\ell ])$ is a subquotient of $
\Hom (M(X)(a), M(Y)(b)[\ell +d_X-d_Y ])= \Hom (M(X\times Y), \Bbb Q (b-a +d_Y )[\ell +d_X +d_Y ])$, which is 0 for $\ell > (2b+d_Y)-(2a+d_X)$ due to the next lemma:

 {\bf Lemma.} {\it  If $X$ is any smooth variety, then $\Hom (M(X), \Bbb Q (n)[\ell])=0$  for  $\ell >n+\text{min}\{ d_X ,n\}$.}

{\it Proof.}  $R\Hom (M(X), \Bbb Q (n))$ is Bloch's complex of relative cycles (see Lecture 19 from \cite{MVW}). Thus $\Hom (M(X), \Bbb Q (n)[\ell ])$ is a subquotient of the group of codimension $n$ cycles on $X\times \Bbb A^{2n-\ell}$, which is 0 for $\ell >d_X +n$ or $\ell >2n$.  \hfill$\square$

{\bf 1.8.}  Suppose an irreducible $P\in \CM_k$ is effective, i.e., occurs in some ${}^\mu\! H^i M(Y)$. By the argument from \cite{Del1}, it occurs then in ${}^\mu\! H^i M(X)$ with $X$ smooth and projective of dimension $\le d_Y$. The {\it level} of $P$ is the smallest dimension of such an $Y$. For any effective $P\in \CM_k$ its level is the maximal level of its irreducible subquotients.

{\bf Proposition.} {\it If $P$, $Q$ are effective  of level $\le \ell$, then $\Hom (P,Q[\ell ])=0$.} 

{\it Proof.} It suffice to check this when $P$, $Q$ are irreducible. Then $P$ occurs  in some ${}^\mu\!H^{-d_X}M(X)(a)$ where $X$ is smooth projective with $d_X\le \ell$ and $0\le a\le \ell-d_X$; same for $Q$. As in the proof in 7,  one can realize $\Hom (P,Q[\ell ])$ as a subquotient of
$\Hom (M(X\times Y), \Bbb Q (b-a +d_Y )[\ell +d_X +d_Y ])$. Now
use the lemma in 1.7.
\hfill$\square$

\bigskip

\centerline{\bf \S 2. Suslin's Lawson homology conjecture yields standard conjectures}

{\bf 2.1.} For a complex projective variety $X$  we have its Lawson homology groups $L_r H_{2r+i }(X,\Bbb Z ):= \pi_i (C_r (X)^+ )$; here $C_r (X)$ is the topological Chow monoid of effective $r$-cycles on $X$, and $C_r (X)^+$ is its group completion. They form the ``homology" part  of a Bloch-Ogus style cohomology theory for complex algebraic varieties, see \cite{Fr}.
There is another cohomology theory with cohomology groups $H^i_\tau (X,\Bbb Z (n)):= H^i (X_{\text{Zar}}, \tau_{\le n} R\pi_* \Bbb Z (n))$; here $\pi : X_{\text{cl}}\to X_{\text{Zar}}$ is the map from the classical topology of $X$ to the Zariski one, $\Bbb Z (n)=(2\pi i )^n \Bbb Z$ is the constant sheaf on $X_{\text{cl}}$, $\tau_{\le n}$ is the truncation. There is a natural morphism from the former cohomology theory to the latter one, and the Suslin conjecture asserts that this morphism is an isomorphism. More concretely, this means that for smooth projective $X$ the canonical map $L_r H_{a }(X, \Bbb Z )\to H_a (X,\Bbb Z (-r))$ is an isomorphism for $a \ge \dim X + r$. 
 
{\it Remark.} The Suslin conjecture with finite coefficients ($\Bbb Z$ replaced by $\Bbb Z/\ell $) is known to be true: indeed, by \cite{SV},   $L_r H_{a }(X, \Bbb Z/\ell  ) $ equals the motivic homology with coefficients in $\Bbb Z/\ell (-r)$,  so the assertion comes from the Milnor-Bloch-Kato conjecture established by Voevodsky, Rost,... Therefore the Suslin conjecture with $\Bbb Z$-coefficients amounts to the assertion that the groups  $L_r H_{a}(X, \Bbb  Z )$ for $a \ge   \dim X + r$  are finitely generated. And
 the Suslin conjecture with $\Bbb Z$-coefficients is equivalent to that with $\Bbb Q$-coefficients.

{\bf 2.2.}  From now on all the (co)homology  have $\Bbb Q$-coefficients, which are omitted,  as well as the Tate twist, in  the notation.

{\bf Proposition.}  {\it The next conjectures are equivalent:\footnote{The implication (iv) $\Rightarrow$ (ii) was observed independently by S.~Bloch and B.~Kahn.} \newline
(i) For any smooth projective $X$ the maps $L_r H_{a }(X)\to H_{a} (X)$ are surjective for $a \ge \dim X + r$. \newline (ii) For $X$ as in (i) one can find a  finite correspondence $f:X\to Y$
with $Y$ projective smooth and $\dim Y =\dim X -1$ such that                 
$f^*: H^{i} (Y) \to H^{i} (X)$ is surjective for $i<\dim X$.
\newline 
(iii) For $X$ as in (i) and  any $j\ge 0$ one can find a  finite correspondence $f_j :X\to Y_j$ with $Y_j$ projective smooth of dimension $j$ such that                 
$f_j^*: H^j (Y_j ) \to H^j (X)$ is surjective. \newline
(iv) The standard conjectures (for varieties over $\Bbb C$).\footnote{Hence over any field of characteristic 0.} }

{\it Proof.}  (iii)$\Rightarrow$(ii): Take $Y=\mathop\sqcup\limits_{j<\dim X} Y_j$, $f=\mathop\Sigma\limits_{j<\dim X} f_j$. (ii)$\Rightarrow$(iii): Take $(Y_{\dim X},f_{\dim X})=(X,\id_X )$. For $j<\dim X$, find
$(Y_j , f_j )$ using downward induction by $j$: namely, $(Y_j , f_j ) = (Y, f f_{j+1})$ where $(Y,f)$ comes from (ii) for $X$ replaced by $Y_{j+1}$. For
 $j> \dim X$ the assertion is evident (say, take $Y_j = X\times \Bbb P^{j-\dim X}$).

(i)$\Rightarrow$(iii): We can assume that  $r= \dim X - j \ge 1$. Recall that $C_r (X)$ is disjoint union of projective varieties, and we have a universal family of $r$-cycles on $X$ parametrized by $C_r (X)$. Viewed as a
correspondence between $C_r (X)$ and $X$, it yields
a map $H_a (C_r (X)) \to H_{a+2r}(X )$, and (i) says that this map is surjective for $a\ge j$. Replace $C_r (X)$ by  its sufficiently large component so that surjectivity still holds. Let $Z$ be a resolution of singularities of the latter; the map $H_a (Z) \to H_{a+2r}(X )$ for $a\ge j$ is still surjective by a usual mixed Hodge theory argument. Let $Y$ be a generic iterated hyperplane section of dimension $j$, so $H_j (Y) \to H_j (Z)$ is surjective by weak Lefschetz, hence $H_j (Y) \twoheadrightarrow H_{j+2r}(X )$, or, replacing the homology by cohomology, $H^j (Y) \twoheadrightarrow H^{j}(X )$. By construction, this map is the action of a correspondence given by a cycle of dimension $\dim X$ on $X\times Y$. By \cite{FV} 7.1, it can be replaced by a finite correspondence,  and we are done.

(iv)$\Rightarrow$(i): Consider  the ``inverse Lefschetz" endomorphism $\Lambda :H^\cdot (X) \to H^{\cdot -2}(X)$.  The map $\Lambda^r : H^{i+2r}(X)\to H^i (X)$ is surjective for $i\le \dim X -r$, i.e., the corresponding map on homology $H_{a-2r}(X)\to H_a (X)$ is surjective for $a\ge \dim X +r$. Realizing $\Lambda^r$ as an $X$-family of $r$-cycles on $X$ (by (iv) and \cite{FV} 7.1), we factor the latter map through $L_r H_a (X)\to H_a (X)$, which yields (i).

{\bf 2.3.}  It remains to prove that (iii) implies (iv). Recall that the standard conjectures reduce to the Lefschetz type conjecture.
                                     
 For a smooth projective variety $X$, consider the next three conjectures
 $L(X)$, $l(X)$, $S(X)$   about the cohomology of $X$:    \newline                                          
- $L(X)$ is the Lefschetz type standard conjecture for $X$; \newline                        
- $l(X)$ is the next assertion: for every $i > 0$ one can find a correspondence on     
$X$ that yields an isomorphism $H^{\dim X + i} (X) \simeq H^{\dim X -i}(X)$; \newline            
- $ S(X)$ is  conjecture (iii) from 2.2 for our $X$ and any $j\ge 0$.

Let $L(n)$ be the assertion that $L(X)$ is true for all $X$ of dimension $\le n$;     
same for $l(n)$, $S(n)$.  We will show that $S(n)$ implies $L(n)$. This takes two steps:  \newline                                                     
(a) $S(n)\,\&\,L(n-1)$ implies $l(n)$, and  (b) $l(X)$ implies $L(X)$.                            
                                                                                 
{\it Proof of} (a).                                                                    
By $S(n)$, we can find smooth projective $Y$ of dimension $n-i$ and a                  
 correspondence $f: X \to Y$ such that $f^* : H^{n-i}(Y) \to H^{n-i}(X)$             
 is surjective. Pick an ample line bundle on $Y$ and consider the                  
 corresponding primitive decomposition of  $H^\cdot (Y)$. By $L(n-1) $               
 the projectors $\pi_a$ on its components are given by algebraic                   
 correspondences. Denote by $ \pi_+ , \pi_-$  the sum of  $\pi_a$'s                     
 such that $\pi_+ + \pi_-$  is the projector onto $H^{n-i}(Y)$ and the                
 Lefschets pairing on the images of $ \pi_+ , \pi_-$ provides a                     
 positively, resp. negatively, defined polarization.                             
                                                                                 
Set $f_\pm := \pi_\pm f$. Thus $f^* : H^{n-i}(Y) \to H^{n-i}(X) $                    
 equals $f_+^* + f_-^*$. Consider the maps                                         
$ f_{\pm *} : H^{n+i}(X) \to H^{n-i}(Y)$.                                          
                                                                                 
{\bf Lemma}. {\it For almost all non-zero rational numbers a the restriction                
 of the Poincar\'e bilinear form to the image of                                  
$ a f_{+*}+  f_{-*}:  H^{n+i}(X) \to H^{n-i}(Y)$ is non-degenerate.  }              
                                                                                 
The lemma implies (a): Indeed, pick $a$ as above; set $f':= af_+ + f_-$.                
Then $f^{\prime *} : H^{n-i}(Y) \to H^{n-i}(X)$ is surjective (since such is               
$f^*$ and $a\neq 0$), hence its adjoint (with respect to the Poincar\'e               
pairings) $f'_* = af_{+*} + f_{-*}: H^{n+i}(X) \to H^{n-i}(Y)$ is                  
injective. The condition of the lemma implies then that                          
$f^{\prime *}f'_*: H^{n+i}(X) \to H^{n-i}(X)$ is an isomorphism, q.e.d.                    

{\it Proof of Lemma}. Consider our cohomology groups with real                         
coefficients. Our picture decomposes into the direct sum of                      
$\Bbb R$-Hodge structure isotypical pieces. It suffices to prove the                    
lemma for one such piece. Our Hodge structures look as                           
$V\otimes H$, where $H$  is a fixed irreducible Hodge structure                       
(rank 2 or rank 1) and $V$ is a real vector space (i.e., a                         
Hodge structure of type (0,0)). If our Hodge structure is a subspace    
of $H^{n-i}(Y)$, then the Poincar\'e pairing is the tensor product  of                
a symmetric bilinear form $q$ on $V$ and a fixed polarization on $H$;                  
if we live in $\pi_\pm H^{n-i}(Y)$, then $q$  is either positive or                   
negative definite. Now the lemma follows  from the next   linear algebra                       
assertion:\footnote{Which follows from the fact that  $U$ can be decomposed into                            
a direct sum of 1-dimensional subspaces  orthogonal                        
with respect to both bilinear forms $g_+^*(q_+ )$  and                              
 $g_-^*(q_- )$.}                                                   
 Let $V_+, V_-$ be $\Bbb R$-vector spaces equipped with, respectively,                    
positive and negative definite symmetric bilinear forms                          
$q_+ $, $q_-$,  and $g_\pm : U \to V_\pm         $                                      
be linear maps; then for almost all non-zero real a the form                     
$q_+ \oplus q_-$ is non-degenerate on the image of                                 
$ag_+ \oplus  g_- : U \to V_+ \oplus V_-$.

{\it Proof of} (b). Assuming $l(X)$, we want to find for every $i>0$ a                     
correspondence $ c$ on $X$  (here $n:=\dim X$)                                           
whose action on $H^\cdot (X)$ is the inverse to Lefschetz                          
$H^{n+i}(X)\to H^{n-i}(X)$, all other components are 0.                             
We do downward induction by $i$. By the induction assumption,                      
all the projectors $ p_j$ on $H^{n+j}(X)$, $|j|>i$, come from                           
correspondences. By $l(X)$, we can find a correspondence $c' $                       
that provides an isomorphism $H^{n+i}(X)\to H^{n-i}(X)$.                            
Multiplying $c'$ by the product of $(1- p_j )$, $j>i$, from the right                   
and by the product of $(1- p_j)$, $j<-i$, from the left, we can                      
assume that the isomorphism $H^{n+i}(X)\to H^{n-i}(X)$                              
is the only non-zero component of the action of $c'          $                     
on $H^\cdot  (X)$.                                                                      
                                                                                 
The composition $A$ of $c'$ with the $i$th power of                          
Lefschetz  acts as an automorphism on $H^{n-i}(X)$, and all                     
its other components are 0. Thus there is a polynomial $f$ in $\Bbb Q [t]$ such                       
that $f(A)A$ acts as identity on $H^{n-i}(X)$. The promised $c  $                      
is $f(A)c'$.    \hfill$\square$

\end